\documentclass[12pt]{amsart}

\usepackage{latexsym,amscd,amsmath,verbatim,amssymb,calc,enumitem,amsxtra,multicol,
mathrsfs,amsthm,ifpdf,xcolor,xspace,mathtools}
\usepackage[all]{xy}
\usepackage{comment}
\usepackage{graphicx, relsize}   
\usepackage[margin=1.50in]{geometry}
\ifpdf
  \usepackage[breaklinks=true]{hyperref}
\else
  \usepackage[hypertex]{hyperref}
\fi

\usepackage{scalerel,stackengine}

\theoremstyle{plain}
\newtheorem{theorem}{Theorem}[section]
\newtheorem{corollary}[theorem]{Corollary}

\newtheorem{lemma}[theorem]{Lemma}
\newtheorem*{theorem*}{Theorem}
\newtheorem*{theoremA*}{Theorem A}
\newtheorem*{theoremB*}{Theorem B}
\newtheorem*{theoremC*}{Theorem C}
\newtheorem*{theoremD*}{Theorem D}
\newtheorem*{theoremE*}{Theorem E}
\newtheorem*{proposition*}{Proposition}

\theoremstyle{remark}
\newtheorem{remark}[theorem]{Remark}

\newtheorem{notation}[theorem]{Notation}

\newtheorem{discussion}[theorem]{Discussion}

\theoremstyle{definition}
\newtheorem{definition}[theorem]{Definition}
\newtheorem{example}[theorem]{Example}

\numberwithin{equation}{theorem}

\newcommand{\new}{\newcommand}
\new{\comm}{\iffalse}
\renewcommand{\comm}{\iftrue}
\new{\op}{\operatorname}
\new{\0}{^{\circ}}
\new{\Ann}{\mathrm{Ann}}
\new{\arwf}{\xrightarrow}
\new{\asb}{asymptotically sparse boundary\xsp}
\new{\Bb}{{\mathbb{B}}}
\new\bem[1]{\bibitem[#1]{#1}}
\new{\ben}{\begin{enumerate}} \new{\een}{\end{enumerate}}
\new{\bena}{\ben[label=$({\rm \alph*})$]}
\new{\benr}{\ben[resume]}
\new\bensp[1]{\ben[label=#1]}
\new{\blank}{\underline{\hbox{\ \ }}}
\new{\card}{\op{card}} 
\new{\cF}{\mathcal{F}}
\new{\CM}{Cohen-Macaulay\xspace}
\new{\cN}{\mathcal{N}}
\new{\congsm}{\scalebox{.8}[.8]{$\cong$}}
\new{\content}{\mathfrak{h}}
\new{\e}{\op{e}}
\new{\fA}{\mathfrak{A}}
\new{\fB}{\mathfrak{B}}
\new{\fb}{\mathfrak{b}}
\new{\fE}{\mathfrak{E}}
\new{\fg}{finitely generated\xspace}
\new{\fG}{\mathfrak{G}}
\new{\fH}{\mathfrak{H}}
\new{\fm}{\mathfrak{m}}   
\new{\fM}{\mathfrak{M}}
\new{\fp}{\mathfrak p}
\renewcommand{\geq}{\geqslant} \renewcommand{\leq}{\leqslant} 
\new{\HH}{\op{H}}
\new{\hl}{{\op{h}}}
\new{\imp}{\Rightarrow}
\new{\inc}{\subseteq} 
\new{\inj}{\hookrightarrow}
\new{\iso}{\overset{\raisebox{-6pt}{\congsm}}{\longrightarrow}  }
\new{\It}{I^{[t]}}
\new{\Ker}{\op{Ker}}
\new{\la}{\lambda}
 \new{\La}{\Lambda}
 \new{\m}{{\mathfrak{m}}}
 \new{\Max}{\op{Max}}
\new{\Min}{\op{Min}}
\new{\Mod}{\op{Mod}}
\new{\n}{\mathfrak{n}}
\new{\N}{{\mathbb{N}}}
\new{\noi}{\noindent}
\new{\ov}{\overline}
\new{\q}{\mathfrak{q}}
\new{\Q}{\mathbb{Q}}
\new{\qL}{\mathcal{L}}
 \new{\Rad}{\op{Rad}}
\new{\rank}{\op{rank}}
\new{\resp}{respectively\xspace}
\new{\sF}{\mathscr{F}}
\new{\sG}{\mathscr{G}}
\new{\sop}{system of parameters\xspace}
\new{\sm}{\,\mbox{\raisebox{1.2pt}{$\smallsetminus$}}\,}
\new{\tM}{\wt{M}}
\new{\tR}{{\widetilde{R}}}
\new{\und}{\underline}
\new{\ua}{{\und{a}}}
\new{\ub}{{\und{b}}}
\new{\ut}{{\und{t}}}
\new{\ux}{\underline{x}}
\new{\uX}{\und{X}}
\new{\uy}{\underline{y}}
\new{\uz}{\underline{z}}
\new{\vect}[2]{{#1}_1, \, \ldots, \, {#1}_{#2}}
\new{\vitem}[2]{%
  \item[]\noindent
  \makebox[\linewidth]{%
    \makebox[0pt][l]{$\vcenter{\hbox{#1}}$}%
    \hfil
    $\displaystyle\begin{aligned}#2\end{aligned}$%
    \hfil
  }
}  
\new{\vv}{\und{1}}
\new{\wh}{\widehat}
\new{\wrt}{with respect to\xsp}
\new{\wt}{\widetilde}
\new{\xsp}{\xspace}
\new{\Z}{\mathbb{Z}}
\new\cred[1]{{\color{red}#1}}

\definecolor{ltgr}{RGB}{211,211,211}  
\new{\blue}[1]{\xspace{\color{blue} #1}\xspace}
\new{\orange}[1]{{\color{orange} #1}}
\definecolor{blgr}{RGB}{31,80,99}
\definecolor{grn}{rgb}{0.0, 0.5, 0.0}
\definecolor{gre}{rgb}{0.0, 0.6, 0.0}
\definecolor{ltgr}{RGB}{150,150,150}  
\definecolor{prp}{RGB}{138,43,226}
\definecolor{ora}{rgb}{0.8, 0.4, 0.0}
\new{\gre}[1]{{\color{gre} #1}}
\new{\mgn}[1]{{\color{magenta} #1}}
\allowdisplaybreaks[1]

\title{Parameters and lim CM sequences:\\ the Q-sequence conjecture and\\ Lech's formula}

\author{Mel Hochster and Olivia Strahan}

\date{\today} 

\thanks{Hochster and Strahan were partially supported by National Science Foundation grant DMS--2200501, and Strahan was also supported by DMS--2101075.}

\subjclass[2020]{Primary. 13H15, 13D22, 13D45, 13D07, 13C14, 13A35}

\keywords{Lech's theorem, lim Cohen-Macaulay sequence, local ring, Q-sequence, system of parameters}

\begin{document}

\begin{abstract} Let $\underline{x}:=x_1,\ldots,x_d \in R$, a ring,
 let $I=(\underline{x})R$, and let $\It=(x_1^t,\,\ldots,\,x_d^t)R$. Let
$\lambda(M)$ denote the length of $M$ over $R$.  $\underline{x} \in R$ is a Q-sequence  
if for all $t \geq 1$, every filtration of $R/\It$ with factors that are quotients of $R/I$ has at 
least $t^d$ factors. Suppose from now on that $\underline{x}$ is a system of parameters for
a local ring $(R,\, \mathfrak m)$.

We prove that every system of parameters $\underline{x}$ of $R$ is a Q-sequence whenever 
$R$ has a weakly lim Cohen--Macaulay sequence and, hence, whenever it has a 
lim Cohen--Macaulay sequence.

Monomial ideals in  $\underline{x}$ are given as
extensions $\fA R$ of monomial ideals $\mathfrak{A} \inc A := \Z[X_1,\, \ldots, \, X_d]$ 
under the map $A \to R$ with $X_1, \, \ldots, \, X_d\mapsto \underline{x}$. 
Let ${\mathfrak B}(\fA)$  denote the set of monomials in 
$A \,\mbox{\raisebox{1.2pt}{$\smallsetminus$}}\, \mathfrak A$ and 
let $\mathcal{N}(\mathfrak{A})$
be  $\#\big({\mathfrak B}({\fA})\big)$, where 
$\#$ is cardinality. Let  $\mathcal{N}_0(\fA)$ be the number of monomials 
in ${\mathfrak B}({\fA})$ in which at least one $X_i$ has exponent 0.

We also prove an extension of Lech's formula: 
with $\underline{x}$ and $I := (\underline x)R \inc  R$ 
as above, if $I$ has multiplicity $e= \operatorname{e}_I(R)$, 
then for every  infinite sequence ${\fA}_t$  of 
$(\underline{X})$-primary  monomial ideals of $A$ such that 
$\gamma(t) := \cN_0({\mathfrak{A}}_t)/\mathcal{N}({\mathfrak{A}}_t) \to 0$ as 
$t \to \infty$, we have that  
$\lim_{t \to \infty} \lambda(R/{\fA}_tR)/\cN({\fA}_t) = e$. 
Examples are given that show that the condition on $\gamma(t)$ cannot be omitted.
\end{abstract}

\maketitle

\pagestyle{myheadings}  
\markboth{MEL HOCHSTER AND OLIVIA STRAHAN}{PARAMETERS, LIM CM SEQUENCES, Q-SEQUENCES, AND LECH'S FORMULA}

\section{Introduction}

Throughout this paper, all rings are commutative, with multiplicative identity.
By a {\it local ring } $(R,\, \m, \, K)$  we mean a Noetherian ring $R$ with unique maximal 
ideal $\m$ and residue class field  $K :=  R/\m$.  

Lim Cohen-Macaulay sequences of modules, which we refer to more briefly 
as {\it lim CM sequences}, were introduced in \cite{BHM26}. See also \cite{Ho17, Ma17, Ma23}.   
The notion of a Q-sequence was introduced in \cite{HoHu09}
and studied further in \cite{HoZh18, Ji21, Ji25, Str25}, and \cite{HoJi26}.  
In \cite{HoHu09} it was proved that 
in equal characteristic, every system of parameters is a Q-sequence,  using 
reduction to characteristic $p$, and it was conjectured that this is true
for every system of parameters in mixed characteristic as well.  This is 
still an open question, even in Krull dimension 3.  This conjecture immediately 
implies the direct summand theorem 
\cite{Ho73, Ho83, And18, Bha18}  and is true for any local ring 
that has a small Cohen-Macaulay module. The existence of big 
Cohen-Macaulay algebras\footnote{See \cite{Ho94, HoHu92, HoHu95, Ho02} for 
equicharacteristic developments and the mixed characteristic case in dimension 3, 
and  \cite{And20, Bha20} for the  mixed characteristic result in all dimensions.}
does not appear to be enough, even in the case of a complete local domain of Krull dimension 3. 

It is proved in this paper that if a local ring $R$ has a lim CM sequence of modules, then 
every \sop is a Q-sequence. It suffices that $\widehat R/P$ has such a sequence for a minimal prime $P$ with $\dim(\widehat R/P)=\dim R$. More generally, it suffices that there exist a local map
of local rings $(R,\, \m) \to (D,\, \n)$ such that $\dim(D) = \dim(R)$, the closed fiber
$D/\m D$ is zero-dimensional, and $D$ has a lim CM sequence. It also suffices if the sequences
of modules are weakly lim CM.

See Theorem~\ref{weaklimQ} and Corollary~\ref{transferQ}.
The conjecture that every complete local domain with algebraically closed 
residue class field has a lim Cohen-Macaulay sequence is, {\it a priori}, much weaker 
than the existence of small \CM modules for these rings.
For example, it is proved in \cite{BHM26} that
lim CM sequences exist in prime characteristic $p>0$ if the local ring 
is F-finite.

We also prove a uniform generalization of Lech's formula under
a mild condition on what might be called the {\it boundary} monomials 
in the ideals in the sequence used. 
See Theorem~\ref{genLECH}. The same theorem gives a zero limit 
when the number of generators exceeds the module dimension. For the 
classical formula, see \cite[Thm.\,11.2.10,\,p.\,229]{SwHu06} and \cite{Lech57}. 

\section*{Acknowledgment}

We would like to thank Craig Huneke for discussion of the results generalizing Lech's theorem.

\section{Background and notation}

We write $\N_+ \inc \N \inc \Z \inc \Q \inc \mathbb{R} \inc \mathbb{C}$ for the positive integers, the nonnegative integers, the integers, 
the rational numbers, the real numbers and the complex numbers, \resp.  All the monomials that we
discuss here have only exponents that are in $\N$.

\begin{definition}\label{kz} Let $\la_R(M) := \la(M)$ be the length of 
the $R$-module $M$ and $\nu_R(M) := \nu(M)$ be the 
least number of 
generators of $M$.  We denote the Koszul homology of 
$M$ \wrt a 
sequence of  elements  $\ux := \vect x{s}$ in $R$ by 
$\HH_{\bullet}(\ux;\,M)$ and we let $\chi^R_i(\ux; \, M) := \chi_i(\ux;\, M)$  denote 
$\sum_{j=i}^{s} (-1)^{j-i}\hl_j(\ux;M)$, where
$\hl_i(\ux;M):=\la_R(\HH_i(\ux;M))$, whenever the lengths in the summation 
are finite. If $i = 0$, the subscript 0 may be omitted, i.e., 
$\chi^R(\ux; \, M) := \chi(\ux; \, M) := \chi^R_0(\ux; \, M)$.  If $i \geq 1$,
$\chi_i(\ux;\,M)$ is called a {\it truncated} Euler characteristic. \end{definition} 

There is a comprehensive list of properties of Koszul homology in \cite[Thm.\,3.4]{BHM26}.

\smallskip


\begin{definition}\label{mult} If  $(R, \, \m, \, K)$ is local, $M$ is a 
nonzero \fg $R$-module with $\dim M=d$, and $I \inc \m$ is an ideal with $s$ generators 
such that $\la(M/IM)<\infty$, where,
necessarily, $s \geq d$, we write
\[ \e_I^s(M) := s! \lim_{t \to \infty} \frac{\la(M/I^tM)}{t^s}.\] 
For the zero module we set $\e_I^s(0)=0$. Equivalently, for nonzero $M$ 
one can compute the displayed limit over $R/\Ann_R(M)$, where the image 
of $I$ is primary to the maximal ideal. Independently of the number 
of generators of $I$, the usual 
Hilbert-Samuel multiplicity is
\[
\e_I(M):=d!\lim_{t\to\infty}\frac{\la(M/I^tM)}{t^d}
\quad(M\ne0),
\]
and we set $\e_I(0)=0$. Thus, when $s=d$, we have
$\e_I^s(M)=\e_I(M)$. Also, let   $\e(M) := \e_{\m}(M)$, 
and if $M = R$ as well, we may simply write $\e$  for $\e(R)$ 
if it is clear from context what $R$ is.
\end{definition}

We recall the following results of Serre \cite{Ser75}:  see also the  summary result on
properties of Koszul homology \cite[Thm.\,3.4(n)]{BHM26}.
For part (b), see \cite[Thm.\,3.4(o)]{BHM26} and 
\cite{Ser75}, Ch.\,\,IV,\,\,Appendice\,\,II, Thm. on p.\,120 and Exemple,\,\,p.\,\,123.
\smallskip

\begin{theorem}[Serre]\label{serre} Let $(R, \m)$ be local, let $\ux$ consist of $s$ elements that
generate an $\m$-primary ideal $I$ 
and let $M$ be a finitely generated module with $\dim M = d$.

\bena
\item The multiplicity $\e^s_I(M) = \chi(\ux;\,M)$. If $\dim(M) < s$, 
then  $\e^s_I(M) = \chi(\ux;\,M) = 0$.
In particular, if $\ux$ is a system of parameters for 
$M$, i.e., if  $s = d$, then
$\e_I(M) = \chi(\ux;\,M)$. 

\item {\rm (Serre positivity.)} Assume that $s =d$, i.e.,  that $\vect xd$ is a 
system of parameters for $M$, and that
$1 \leq i \leq d$.  Then the truncated Euler characteristic 
$\chi_i(\ux;\, M) \geq 0$, and  $\chi_i(\ux;\,M) > 0$ unless 
$\hl_j(\ux;\,M) =0$ for all $j \geq i$. 
\een \end{theorem}

\smallskip

\begin{discussion}\label{formal}{\bf Formal descriptions of ideals.} It will be convenient to have a way of giving  formal descriptions of ideals generated by monomials in elements $\vect rd$ in a commutative ring $R$.  
Let $A:= A_d := \Z[\vect Xd]$ be a polynomial ring. 

If $\ua := (\vect a d) \in \N^d$, we write
$\uX:= \vect Xd$ and $\uX^\ua := X_1^{a_1} \cdots X_d^{a_d}$ 
for a typical monomial in 
the elements $\uX$. Define 
$\fM_d := \{\uX^\ua : \ua \in \N^d\}$.  We call an ideal
 $\fA \inc A$ generated by a subset
of $\fM_d$ a {\it formal monomial} ideal, which will be helpful in making clear which 
monomial ideals are in $\Z[\uX]$ and which are in $R$.

Given an $(\uX)$-primary monomial ideal $\fA$, let 
$\fE(\fA) := \{\ub \in \N^d: \uX^\ub \not \in \fA\}$, and let 
$\fB(\fA) = \{\uX^{\ub}:\ub \in \fE(\fA)\}$.  Let 
$\cN(\fA) := \#\big(\fB(\fA)\big) = \#\big(\fE(\fA)\big)$.
The residue classes of  the monomials in $\fB(\fA)$ 
form a  free $\Z$-basis for $A/\fA$.

For $\ua, \, \ub \in \N^d$,  $\ua \leq \ub$ 
means $a_i \leq b_i$ for each $1 \leq i \leq d$.
With this notation, \smallskip

\hspace*{30pt}$(\sharp)$\ \ if $\ua, \, \ub \in \N^d$ with 
$\ua \leq \ub$ and $\ub \in \fE(\fA)$, \smallskip
then $\ua \in \fE(\fA)$. 

Let $\fG$ be any 
set of monomial generators of a formal monomial 
ideal $\fA \inc \Z[\uX]$.  It is routine to check that 
$\uX^{\ub} \in \fA$ if and only if $\ua \leq \ub$ for some  $\uX^{\ua} \in \fG$. 
\end{discussion}

\begin{remark} It is elementary to prove that if $\ux$ is a 
system of parameters for a local ring, then the map $\Z[\uX] \to R$ 
induces a bijection from monomials in $\uX$ to monomials
in $\ux$.  It suffices to show this after replacing $R$ by
$\wh{R}/\fp$, where $\fp$ is a minimal prime of $\wh{R}$
with $\dim \wh{R}/\fp = \dim \wh{R}$. By factoring out any formal
common factors, we then obtain an equality 
between the  images of two formally distinct  monomials $\mu = \mu'$ in disjoint sets 
of parameters in a catenary local domain,  and their common value cannot 
be the identity 1. Suppose that $x_j$ occurs in $\mu'$. 
Choose a minimal prime $\q$ of $x_j$. Since $\mu = \mu' \in (x_j) \inc \q$,
one of the parameters $x_i$ with $i \neq j$ occurring in the product $\mu$ 
must be in the prime $\q$.  Since $R$ is catenary, we cannot have two elements of a \sop
in a prime ideal $\q$ of height one. \smallskip

However, it is not an elementary fact that the extension of 
a formal monomial ideal $\fA \inc \Z[\uX]$ contains only those 
monomials in $\ux$ that are in the image of the set of monomials in $\fA$.  
This statement is true, but it is equivalent to the direct summand 
theorem, and we do not want to assume the direct summand
theorem in our treatment here.
\end{remark}

\subsection{Lim Cohen-Macaulay sequences}\label{limcm}
For a local ring $(R,\m,K)$ of Krull dimension $d$, a sequence of
finitely generated modules $M_n$ of dimension $d$ is
\emph{lim Cohen-Macaulay} or \emph{lim CM} if
$\hl_i(\ux;M_n)=\op{o}(\nu(M_n))$ for all $i\geq1$ and one
(equivalently, every) system of parameters $\ux$. 
Over a domain,
$\nu(M_n)$ may be replaced by $\op{rank}_R(M_n)$
\cite[Theorems 4.6 and 4.17]{BHM26}. Every $F$-finite local ring of
characteristic $p>0$ has the lim CM sequence $\{\cF^n_*(R)\}_n$,
where $\cF^n_*$ denotes restriction of scalars along
$r\mapsto r^{p^n}$ \cite[Theorem 5.4]{BHM26}. These sequences were
introduced partly to approach Serre's positivity conjecture:
their existence for every complete mixed characteristic local
domain with algebraically closed residue 
field would settle the remaining ramified case 
\cite[\S1.6\,and\,Thm.\,6.1]{BHM26}. We note that on October 6, 2026, OpenAI 
announced a proof \cite{OpenAI26}, which references \cite{BHM26}, 
of the ramified case of the 
Serre positivity conjecture for multiplicities.

A sequence is called \emph{weakly lim CM} if
$\chi_1(\ux;M_n)=\op{o}(\nu(M_n))$ for one (equivalently, every) system of parameters 
$\ux$ \cite[Proposition~2.6]{Ma23}. If, additionally,
$\e(M_n)/\nu(M_n)\to 1$, the sequence of modules is called
\emph{weakly lim Ulrich}.
Linquan Ma proved that the existence of such a sequence over 
a local domain $R$ implies Lech's inequality $\e(R)\leq \e(S)$ 
for every flat local map $R\to S$ \cite[Definitions 2.2 and 2.9, 
Proposition 2.6, and Theorem 2.11]{Ma23}.

\subsection{Quasilength and Q-sequences}\label{QQ}
Let $I$ be a finitely generated ideal of a ring $R$. The
\emph{$I$-quasilength} $\qL_I(M)$ of $M$ is the least number of factors in
a finite filtration of $M$ whose factors are cyclic and annihilated
by $I$ or, equivalently, are quotients of $R/I$. Set
$\qL_I(M)=\infty$ if no such filtration exists. This invariant is
finite precisely when $M$ is finitely generated and some power of
$I$ annihilates $M$. Concatenating filtrations proves that if one has
an exact sequence $0\to  M'\to M \to M''\to 0$ then 
$\qL_I(M)\leq \qL_I(M')+\qL_I(M'')$.

Taking images of a filtration also shows that $I$-quasilength cannot
increase on passing to a quotient or after extension of scalars;
see \cite[Prop.\,1.1]{HoHu09}.

For $\ux=x_1,\ldots,x_d$, put $I=(\ux)R$ and
$\It=(x_1^t,\ldots,x_d^t)R$. This notation is slightly inaccurate, since $\It$
depends on $\ux$, not just $I$, but we only use it when generators $\ux$ for $I$ have
been specified.  Assume $I\ne R$.
The sequence $\ux$ is a \emph{Q-sequence} if $\qL_I(R/\It)=t^d$ for every $t \geq 1$.
The upper bound $\qL_I(R/\It)\leq t^d$ always holds. More generally,
let $\fA \inc \Z[X_1,\ldots,X_d]$ be a $(\uX)$-primary formal monomial ideal
and let $J=\fA R$ be its extension to $R$. Then
$\qL_I(R/J)\leq \cN(\fA)$. 

\smallskip

To see this, linearly order all of the monomials 
$\mu_1,\, \ldots,\, \mu_{\cN(\fA)}$ of
$\fB(\fA)$ so that their total degrees are descending.
Let $J_0 := J$ and $J_i = J_{i-1} + \mu_i R$ for
$1  \leq i \leq \cN(\fA)$. 
Thus, we have  $J = J_0\, \inc \, J_1 \, \inc\, \cdots \, \inc\, J_{\cN(\fA)} = R$. 
Here, consecutive ideals are allowed to be equal. It is then clear that $x_kJ_{i+1} \subseteq J_i$ for all $i,\, k$. 
The resulting filtration of $R/J$ has
$\cN(\fA)$ cyclic factors annihilated by $I$, before 
omitting  any factors that are 0, as one may.

In fact, $\ux$ is a Q-sequence if and only if
\bensp{$(\dagger)$}
\item  $\qL_I(R/\fA R)=\cN(\fA)$ for every $(\uX)$-primary  formal monomial ideal $\fA$.
\een
Note that the formal description is part of the data: $\cN(\fA)$ is not
defined by merely counting distinct elements of $R$.
For the forward implication, choose $t$ such that
$(X_1^t,\ldots,X_d^t) \inc \fA$. Among the $t^d$ monomials
$X_1^{a_1}\cdots X_d^{a_d}$ with $0\leq a_i<t$, exactly
$t^d-\cN(\fA)$ belong to $\fA$. Adjoining their evaluations to
$\It$, again in decreasing total degree, gives $\qL_I(J/\It)\leq t^d-\cN(\fA)$.
Consequently, the exact sequence $0\to J/\It \to R/\It\to R/J \to 0$ gives
\[
 t^d=\qL_I(R/\It)
       \leq t^d-\cN(\fA)+\qL_I(R/\fA R),
\]
which supplies the reverse inequality. Conversely, apply
$(\dagger)$ with $\fA$ equal to $(X_1^t,\ldots,X_d^t)$.
This is the monomial characterization of
\cite[Remark~3.7]{HoHu09}. 

For comparison with the arguments below, we use the normalized
quasilength invariant
\[
 \content^d_{\ux}(R)
      :=\lim_{t\to\infty}\frac{\qL_I(R/\It)}{t^d}
       =\inf_{t\geq1}\frac{\qL_I(R/\It)}{t^d}.
\]
The limit and its equality with the infimum are proved in
\cite[Theorem~2.4]{HoHu09}, where it is referred to as {\it content}. 
Thus, $0\leq\content^d_{\ux}(R)\leq1$,
and $\ux$ is a Q-sequence precisely when
$\content^d_{\ux}(R)=1$. It is intended to be thought of as a partial indicator
of how likely it is that the local cohomology module $H^d_{(\ux)}(R)$ will be nonzero 
and remain so after a base change.

The Q-sequence conjecture asserts that every system of parameters
of a Noetherian local ring is a Q-sequence. It is known in equal
characteristic and in dimensions at most two
\cite[Theorems~4.1, 4.6, and 4.7]{HoHu09}. The remaining general
problem is in mixed characteristic in dimension at least three;
see \cite[Conjecture~8.1]{HoZh18} and \cite{Str25} for further
mixed characteristic cases.

\section{Proof of the Q-sequence conjecture
in the lim CM case}\label{proof limCM}\label{limCM}

\begin{theorem}\label{weaklimQ}
Let $(R,\m,K)$ be a local ring of dimension $d$ admitting a weakly
lim CM sequence. Then every system of parameters of $R$ is a
Q-sequence. In particular, this holds if $R$ admits a lim CM sequence.
\end{theorem}

\begin{proof}
Fix a system of parameters $\ux$, put $I=(\ux)R$, and let
$\{M_n\}_n$ be the given sequence. Write
$\alpha_n=\la(M_n/IM_n)>0$. By Serre's formula,
\[
 \e_I(M_n)=\alpha_n-\chi_1(\ux;M_n).
\]
Since $\alpha_n\geq\nu(M_n)$ and
$\chi_1(\ux;M_n)=\op{o}(\nu(M_n))$, we obtain
\[
                 \lim_{n\to\infty}\frac{\e_I(M_n)}{\alpha_n}=1.
\]
For $t\geq1$, choose an $I$-quasilength filtration of $R/\It$.
Tensor it with $M_n$ and take images in $M_n/\It M_n$.
Right exactness of tensor shows that each successive quotient is a quotient
of $M_n/IM_n$ and so has length at most $\alpha_n$. Consequently,
\[
 \la(M_n/\It M_n)\leq\qL_I(R/\It)\,\alpha_n.
\]
For each fixed $n$, divide by $t^d \alpha_n$ and let $t$ tend to infinity.
The classical Lech formula, which says, in part, that 
\[\lim_{t\to \infty} \frac{\lambda(M/\It M)}{t^d} = \e_I(M)\] 
for any finitely generated module $M$ of dimension $d$ and parameter ideal $I$
\cite[Theorem~11.2.10]{SwHu06}, gives
\[ \frac{\e_I(M_n)}{\alpha_n}\leq\content^d_{\ux}(R)\leq1.\]
Thus, since the leftmost term approaches $1$ as $n$ tends to infinity, 
we obtain that $\content^d_{\ux}(R)=1$,
which is equivalent to the Q-sequence condition by \S\ref{QQ}.
\end{proof}

\begin{corollary}\label{transferQ}
Let $(R,\m)\to(D,\mathfrak n)$ be a local homomorphism of
Noetherian local rings of the same dimension $d$, and assume that
$\sqrt{\m D}=\mathfrak n$. If $D$ admits a weakly lim CM sequence,
then every system of parameters of $R$ is a Q-sequence.
\end{corollary}

\begin{proof}
The images of a system of parameters $\ux$ of $R$ form a system
of parameters of $D$: some power of $\m$ is contained in $(\ux)R$,
so $\sqrt{(\ux)D}=\mathfrak n$. Let $I = (\ux)R$. Extension of scalars takes an
$I$-quasilength filtration of $R/\It$ to a filtration of
$D/\It D$ with no more factors. Theorem~\ref{weaklimQ} therefore gives \vspace*{7pt}

\hspace*{1.4in} {$t^d=\qL_{ID}(D/\It D)\leq\qL_I(R/\It)\leq t^d.$} \end{proof}

\section{A generalization of Lech's formula for multiplicities}\label{lechdiscussion}

Let $(R,\,\m,\,K)$ be a local ring of dimension $d$. Let $\la_R(\blank)$ 
or simply $\la(\blank)$ 
denote length over $R$. 

A well-known theorem of Lech \cite[Theorem~11.2.10]{SwHu06} asserts that if 
$\vect xd$ is a \sop for a local ring $(R, \m)$, 
$I = (\vect xd)R$, $\ut = (\vect t d)$ is a $d$-tuple of positive 
integers, $I^{[\ut]} := (x_1^{t_1},\, \ldots, \, x_d^{t_d})R$, and
$M$ is a \fg  $R$-module, 
then
\[ \lim_{\ut \to \infty} \frac{\lambda(M/I^{[\ut]}M)}{t_1\cdots t_d} = \e_I(M)\]
whenever $\dim(M) =d$, while 
the limit is 0 if $\dim(M) <d$.  
The limit is taken so that all of the
$t_i$ tend to infinity. Note that the special case where the $t_i$ are all the same
was used earlier, in the proof of Theorem~\ref{weaklimQ}.

We consider here a generalization of Lech's theorem, Theorem~\ref{genLECH},
that allows a much larger family of sequences of monomial ideals to play the 
role of the  $\It$ or the $I^t$ in formulas for multiplicities.

\begin{notation} Let $\uX$ be indeterminates over $\Z$. We freely 
use the terminology and notation, especially $\fB(\fA)$, $\fE(\fA)$, and 
$\cN(\fA)$,  introduced in Discussion~\ref{formal}.
In this section, $\#(\uX)$ is typically either $d$ or
$s \geq d$. \smallskip

We refer to monomials $\uX^{\ua}$ as 
{\it boundary} monomials if one or more of the $a_i$ is 0, because the 
elements $\ua$ with one or more entries equal to 0 form the boundary of the first
orthant. We use $\cN_0(\fA)$ to denote the number of boundary monomials in $\fB(\fA)$.
\end{notation}

\begin{discussion}\label{boundary}{\bf A key boundary condition.} 
Fix $s \geq 1$, and let $\uX : = \vect X s$.  We introduce a 
condition on the relative number $\gamma(\fA)$
of boundary monomials in $\fB(\fA)$, defined formally just below,
that plays a critical role in the generalization of Lech's formula. 
\smallskip                                               

Fix an $(\uX)$-primary formal monomial ideal 
$\fA\ \inc \Z[\uX].$  For every $1 \leq i \leq s$, we define 
$\beta_i(\fA) := \#(\{\ua\in \fE(\fA):a_i=0\})$.  
Let $\cN_0'(\fA) := \sum_{i=1}^s \beta_i(\fA)$. It is clear that
\[ (\natural)\quad \cN_0(\fA) \leq \cN_0'(\fA)  
\leq s \max\{\beta_i(\fA): 1 \leq i \leq s\} \leq s\,\cN_0(\fA).\]                           
We define
$\displaystyle  \gamma(\fA) := \frac{\cN_0(\fA)}{\cN(\fA)}$ and
$\displaystyle \gamma'(\fA) := \frac{\cN_0'(\fA)}{\cN(\fA)}$.    

We say that a sequence or directed system of 
$(\uX)$-primary formal monomial ideals $\fA_t \inc \Z[\uX]$, 
for $t$ in a suitable index set, has {\it \asb} if for all 
$\epsilon > 0$ in $\mathbb{R}$, there exists an index $t_{\epsilon}$ 
such that $\gamma(\fA_t) \leq \epsilon$ for
all $t \geq t_{\epsilon}$. We also write this condition as 
$\lim_{t \to \infty} \gamma(\fA_t) = 0$ or as   
$\gamma(\fA_t) \to 0$ as $t \to \infty$.

Evidently, if we replace $\cN_0(\fA_t)$ by any of the four terms in the
display $(\natural)$, with $\fA_t$ replacing $\fA$ throughout, we get 
an equivalent condition, since $(\natural)$ implies that
the ratio of any two of these four terms lies between $1/s$ and $s$. 
In particular, a sequence has asymptotically sparse boundary if and only if
$\gamma'(\fA_t)\to 0$ as $t\to \infty$, which 
is helpful for some specific calculations.
\end{discussion}   

Recall Definitions~\ref{kz} and~\ref{mult} and Theorem~\ref{serre}.  We prove a
preliminary result, and then state and prove our generalization of Lech's 
formula, which is Theorem~\ref{genLECH}.
The next two lemmas do nearly all the work necessary to prove the 
main result Theorem~\ref{genLECH} of this section.

\begin{lemma}\label{pre-main} Let $\uX :=\vect X d$ and $\fA \inc \Z[\uX]$ be 
a $(\uX)$-primary monomial ideal in $\uX.$ Let $M$ be a \fg module
over the local ring $(R,\,\m)$, where $M$ and $R$ both have Krull dimension $d$, 
and let  $\Z[\uX] \to R$ be such that $\uX \mapsto \ux$, where $\ux := \vect xd$
is a system of parameters for $R$. Let $I:=(\ux)R$. 
Then $\la(M/\fA M)\geq \cN(\fA)\e_I(M)$.                 
\end{lemma}
\begin{proof} Let $y_i$ denote the image of $X_i$ in $\Z[\uX]/\fA$. Define 
$\tR := R \otimes_{\Z}(\Z[\uX]/\fA)$.  Then $\tR$ is a free $R$-module of rank
$\cN(\fA)$. In fact, if we abuse the notation slightly by continuing to write $y_i$ for its image $1\otimes y_i$ in $\tilde R$, then $\tilde R$ is generated by $\uy = \vect y d$ as an $R$-algebra and the elements $y_i$ are nilpotent. Moreover, the substitution $\uX \mapsto \uy$ gives a bijective correspondence between $\fB(\fA)$ and a free basis for $\tilde R$ over $R$ that contains $1$.

Thus, $\tR/(\uy) \cong R$ is local, $(\uy)\tR$ is nilpotent,  and $R$ is an 
algebra retract of $\tR$. Both local rings have dimension $d$.  The maximal ideal 
of $\tR$ is $(\m,\, \uy)\tR = \m \oplus_R (\uy)\tR$, and $K$ is the residue field of both
$R$ and $\tR$.

Now let $z_i:=x_i-y_i$.  One of the usual 
monomial filtrations of $\tR$ has $\cN(\fA)$ factors isomorphic 
to $R$ on which every $y_i$ acts as zero. Such a filtration is achieved by making 
a suitable linear ordering of the specified basis of monomials in the $y_i$ that are images
of monomials in $\fB(\fA)$ such that total degree is decreasing and 
using the $R$-free submodules spanned by the initial segments of this ordered free basis.
Consequently, the filtration is 
split as a filtration of $R$-modules, so that tensoring it with $M$ 
gives a filtration of $\tM:= \tR \otimes_R M$ with $\cN(\fA)$ factors 
isomorphic to $M$.  On each factor, $z_i$ acts as $x_i$.

Since $\tM/(z_1,\ldots,z_{d})\tM \cong M/\fA M$
under the substitution $\uy \mapsto \ux$, this quotient has finite length.  
Additivity of Koszul Euler characteristics along the filtration gives that 
$(*) \ \ \chi^{\tR}(\uz;\,\tM) = \cN(\fA)\chi^R(\ux;\,M)= \cN(\fA)\e_I(M)$.
Lengths over $\tR$ and $R$ agree for  finite length $\tR$-modules regarded 
as $R$-modules, since the residue fields are the same. Using $(*)$  
and Serre's nonnegativity for Koszul $\chi_1$, Theorem~\ref{serre}(b), we have that\newline

$ \displaystyle
 \la(M/\fA M)-\cN(\fA)\e_I(M)
 = \hl_0(\uz; \tM)-\chi^{\tR}(\uz;\tM)  =\chi^{\tR}_1(\uz;\,\tM)\geq 0
$.
\end{proof}

\begin{lemma}[Key Lemma]\label{key} Let $\uX :=\vect X s$ 
and let $\fA$ be an $(\uX)$-primary formal monomial ideal 
in $\Z[\uX]$. Let $M$ be a \fg $R$-module, where $(R,\, \m, \, K)$ is a local ring, and
let $\ux := \vect xs \inc \m$, let $I := (\ux)R$, and assume that $M/(\ux)M$ has finite length. View $R$ and $M$ as $\Z[\uX]$-modules by restriction of scalars along $\uX \mapsto \ux$.
As usual, let $I^{[k]} = (x_1^k, \, \ldots, \, x_s^k)R$.
Then for every integer $k \geq 1$ we have that: 
\[(\ddagger)\quad \e_I^s(M)\leq\frac{\la(M/\fA M)}{\cN(\fA)}
 \leq\frac{\la(M/I^{[k]} M)}{k^s} + (k-1)\gamma(\fA)\la(M/IM).\] 
\end{lemma} 

\begin{proof} 
The case $M=0$ is immediate. Otherwise work over
$R/\Ann_R M$. Now $I$ is primary to the maximal
ideal, $\dim R=\dim M=d\leq s$, and all the lengths are
unchanged. If $s>d$, the lower bound is clear,
since $\e_I^s(M)=0$ by Serre's results Theorem~\ref{serre}. 
If $s=d$, use Lemma~\ref{pre-main}.

Let $A :=\Z[\uX]$.  
Fix $k\geq1$.  Partition $\N^s$ 
into $s$-dimensional boxes with sides of length $k$, each  
containing $k^s$ points of $\N^s$, namely:
\[
\Bb^k_{\ua} := \{k\,\ua+\ub: 0 \leq b_i \leq k-1,\ 1 \leq i \leq s\},\quad \ua \in \N^s
\]

Let $\eta_k$ be the number of these boxes completely contained in $\fE(\fA)$, which we
call {\it full} boxes. Let $\rho_k$ be the number of  points of $\fE(\fA)$ in the 
remaining boxes, which we call {\it partial} boxes. Then 
\ben
\item $\displaystyle \eta_k k^s+\rho_k=\cN(\fA) 
\text{\ and\ } \eta_k = \frac{\cN(\fA) - \rho_k}{k^s}$. 
\een
Note that $k\ua \leq \ub$ for any other point $\ub \in \Bb_{\ua}^k$, so we can visualize 
$k\ua$ as the ``lowest" corner of the box $\Bb^k_{\ua}$. 
Linearly order the boxes meeting $\fE(\fA)$ so that the
sums $a_1 + \cdots + a_s$ corresponding to their indices
$\ua$ are descending, breaking ties arbitrarily. Starting with $\fA$,
adjoin, box by box in this order, all monomials whose exponent
vectors lie in $\Bb^k_{\ua}\cap\fE(\fA)$.
The resulting monomial spans are ideals: multiplying a monomial
by a variable either keeps its exponent vector in the same box,
puts it in an earlier box, or puts the monomial in $\fA$.
If we work modulo $\fA$, we obtain a filtration of $A/\fA$.

The factor $\sF_\ua$ corresponding to a full box 
$\Bb^k_{\ua}$ is a cyclic module generated by $\uX^{k\ua}$.
There is an isomorphism $\Z[\uX]/(\uX)^{[k]} \iso \sF_\ua$ under 
the mapping $1 \mapsto \uX^{k\ua}$. Refine the contribution to the filtration from a partial box $\Bb^k_{\ua}$ by adjoining the monomials $\uX^{\ub}$ with
$\ub\in\Bb^k_{\ua}\cap\fE(\fA)$ one at a time, in descending
total degree. Each resulting factor
is isomorphic to $A/(X_1,\ldots,X_s)$. Consequently there are 
$\eta_k$ factors isomorphic to $A/(\uX)^{[k]}$ and $\rho_k$ 
factors isomorphic to $A/(\uX)$.

Regard $R$ and $M$ as $A$-modules via $X_i\mapsto x_i$. Tensor the above 
filtration with $M$ over $A$ and take images in
$(A/\fA)\otimes_A M=M/\fA M$.
The new factors are quotients of the tensor product of $M$ with the corresponding 
original factor, by the right exactness of tensor. This gives a filtration of 
$M/\fA M$ where $\eta_k$ factors are quotients of 
$M \otimes_A (A/(\uX)^{[k]}) \cong M/I^{[k]}M$ 
and $\rho_k$ factors are quotients of $M \otimes_A (A/(\uX)) \cong M/IM$.
Therefore
\benr
\item $\la(M/\fA M)\leq\eta_k\la(M/I^{[k]}M)+\rho_k\la(M/IM)$.  
\een

Note that $\displaystyle \frac{\eta_k\la(M/I^{[k]}M)}{\cN(\fA)}
\leq\frac{\la(M/I^{[k]}M)}{k^s}.$ If $\eta_k=0$, the inequality is obvious, and otherwise it follows from $\eta_k k^s \leq \cN(\fA)$. Consequently, the upper bound in 
$(\ddagger)$ will follow from 
(2) after we bound $\rho_k$.
Let $\vv$ denote the vector $(1,\,1,\, \ldots,\, 1) \in \N^s$.
If $\ub\in \fE(\fA)$ belongs to a
partial box, then $\ub+(k-1)\vv \notin \fE(\fA)$: otherwise
downward closure of $\fE(\fA)$ would put the top corner, hence all of
that box, in $\fE(\fA)$. Hence the partial-box points belong to the set
$\fE(\fA)\sm\big(\fE(\fA)-(k-1)\vv\big)$. Since $\fE(\fA)$ and $\fE(\fA)-(k-1)\vv$ 
have equal finite cardinality, 
\[
 \rho_k\leq\#\big(\fE(\fA)\sm\big(\fE(\fA)-(k-1)\vv\big)\big)
 =\#\big((\fE(\fA)-(k-1)\vv)\sm \fE(\fA)\big).
\]
From the downward closure condition $(\sharp)$ on $\fE(\fA)$ given in 
Discussion~\ref{formal}, we  
have that $\big(\fE(\fA)-(k-1)\vv\big)\cap\N^s\subseteq \fE(\fA)$.
Thus the last cardinality counts exactly those $\ua\in \fE(\fA)$
with $\ua - (k-1)\vv \notin \N^s$, or equivalently, $\min_i a_i <k-1 $. 
For each $0\leq j < k-1$, the map $\ua \mapsto \ua - j\,\vv$ provides an injection
\[\{\ua\in \fE(\fA):\min_i a_i=j\} \inj  \{\ub\in \fE(\fA):\min_i b_i=0\}.\]
$\begin{aligned} {\hbox{\rm It\ follows\ that\qquad } \rho_k}  &\leq \#\big((\fE(\fA)-(k-1)\vv)\sm \fE(\fA)\big) \\
               & =\#\big(\{\ua \in \fE(\fA) : \min_i a_i < k-1\}\big) \\
              &\leq (k-1)\#\big(\{\ub \in \fE(\fA) : \min_i b_i = 0 \}\big) \\     
              & = (k-1) \cN_0(\fA).  \end{aligned}$

Combining this result with $(\natural)$ in Discussion~\ref{boundary}, we have
\benr 
\item $\rho_k \leq (k-1)\cN_0(\fA) \leq (k-1)\sum_{i=1}^s \beta_i(\fA)$ 
\een  

For $k=1$ the same inequality holds because every box is a
singleton and $\rho_1=0$. Combining (1), (2), and (3) proves the upper bound.
\end{proof}

\begin{theorem}\label{genLECH}
Let $(R,\m,K)$ be a Noetherian local ring, let $M$ be a nonzero
finitely generated $R$-module, and let $\ux : = x_1,\ldots,x_s\in \m$.
Let $I=(\ux)R$,  and assume that $\lambda(M/IM)<\infty$. 
As usual, let $\It := (x_1^t, \, \ldots, \, x_s^t)R$, let $\uX := \vect X s$, 
and view $R$ and $M$ as $\Z[\uX]$-modules via $\uX \mapsto \ux$. 
Let  $d=\dim M\leq s$.

For every directed family of $(\uX)$-primary formal monomial ideals 
$\{\fA_t\}_t$ with \asb, we have the generalized Lech formula:
\[ \lim_{t\to\infty}\frac{\la(M/\fA_t M)}{\cN(\fA_t)}
 =\begin{cases}\e_I(M)&\text{if }\dim M=s,\\0&\text{if }\dim M<s.\end{cases}
\]
\end{theorem}
\begin{proof} We may replace $R$ by
$R/\Ann_R(M)$ without changing the lengths or multiplicities, so
that $\dim R=d$ and $I$ is primary to the maximal ideal.
If $s=d$, the usual Lech formula
\cite[Theorem~11.2.10]{SwHu06} gives
$\la(M/I^{[k]}M)/k^s\to\e_I(M)$.
If $s>d$, the containment
$I^{s(k-1)+1}\inc I^{[k]}$ gives
\[0\leq\la(M/I^{[k]}M) \leq \la(M/I^{s(k-1)+1}M)=O(k^d),\]
by the Hilbert-Samuel polynomial, so that division by $k^s$ gives
limit $0$. Thus, in either case, the first summand on the right
of $(\ddagger)$ tends to $\e_I^s(M)$.

Let $\epsilon > 0$  be given.  First fix 
$k \geq 1$ so that 
$\displaystyle \left|\frac{\la(M/I^{[k]} M)}{k^s} - \e^s_I(M) \right| 
< \epsilon/2$. 
Since we are given that $\gamma(\fA_t) \to 0$
as $t \to \infty$, we may choose $t_0$ so that for all $t \geq t_0$,
$(k-1)\gamma(\fA_t)\lambda(M/IM) < \epsilon/2$.

Hence, for every $\epsilon > 0$, there is a 
$t_0$ such that for all $t \geq t_0$, we have that\ \ 
$\displaystyle \left|\frac{\la(M/\fA_t M)}{\cN(\fA_t)}-\e_I^s(M)\right|<\epsilon$. 
Thus, the limit is as stated.
\end{proof}

\begin{corollary}\label{LECH-s-big}
Let $(R,\m)$ be local of dimension $d$, let
$x_1,\ldots,x_s \in \m$ generate an $\m$-primary ideal, and assume
$s>d$. For every sequence of $(\uX)$-primary formal monomial ideals
$\fA _t \inc \Z[X_1,\ldots,X_s]$ with $\gamma(\fA _t)\to 0$ and
every finitely generated $R$-module $M$, one has
\[
 \lim_{t\to\infty}\frac{\la_R(M/\fA _t M)}{\cN(\fA_t)}=0.
\]
In particular, this holds for $M=R$.
\end{corollary}

\begin{proof}
If $M \neq 0$, then $\dim M \leq d < s$, so the theorem applies.
The assertion for $M=0$ is immediate.
\end{proof}

\begin{remark} If $M$ is Cohen--Macaulay and $s=\dim M$, the bounds 
in Lemma~\ref{key} with $k=1$ give the exact equality
$\la(M/\fA M)=\cN(\fA )\e_I(M)$ for every $(\uX)$-primary monomial ideal $\fA$.
No boundary condition is needed in this case. 
\end{remark}

\begin{remark} Let $\uX := \vect X s$ and $\ut := \vect ts \in \N_+^s$.
Consider the family of formal monomial ideals $\fA_{\ut} := (\uX)^{[\ut]} \inc \Z[\uX]$ directed  
by $\leq$ applied to vectors of exponents $\ut$ which are the subject of 
Lech's original formula. One has $\gamma'(\fA_{\ut})=\sum_i 1/t_i$.  
Thus, Theorem~\ref{genLECH} recovers the classical result of Lech as all of the $t_i$ tend 
to infinity.  For the sequence 
$\{\fA_t\}_t$, where $\fA_t := (\uX)^t$ and $t \in \N_+$, 
the corresponding data are:
\[
 \cN(\fA_t)=\binom{t+s-1}{s},\qquad
 \beta_i(\fA_t)=\binom{t+s-2}{s-1},\qquad
 \gamma'(\fA_t)=\frac{s^2}{t+s-1} \to 0.
\]
Therefore, powers of the parameter ideal are also covered by Theorem~\ref{genLECH}.
\end{remark}

\section{The need for a condition on boundary monomials}

The examples in this section show that the \asb condition cannot be omitted in a theorem 
that is to hold for every local ring and module. The first example is quite simple, but it
is not reduced. In the second example the rings are complete local domains.
Throughout this section let $K$ denote an arbitrary field. The examples are valid in
all characteristics.

\begin{example}\label{nonreduced} 
Let $C := K[[x_1,\ldots,x_s]]$.  For all $i$, $1 \leq i \leq s$, consider
$R_i=C[W]/(W^2,x_iW).$ This is a local ring with residue field $K$, and, 
as a $C$-module, $R_i\cong C\oplus W(C/x_iC).$
It has dimension $s$, the elements $x_1,\,\ldots,\,x_s$ form a system of parameters, and 
$\e_{(\ux)}(R_i)=1$.  Moreover, 
$\la_{R_i}(R_i/\fA R_i)=\cN(\fA)+\beta_i(\fA)$.
Thus, if a sequence of $(\uX)$-primary monomial  ideals $\fA_t$ gives the asserted 
multiplicity limit for every local ring, applying it to these $R_i$ forces $\beta_i(\fA_t)/\cN(\fA_t)\to0$ 
for every $i$. Therefore, $\gamma(\fA_t) \to 0$.  
The zero limit assertion in Corollary \ref{LECH-s-big} yields the same necessity 
even more directly: for $R=C/x_iC$ of 
dimension $s-1$, the images of all $s$ variables generate the maximal ideal and 
$\la(R/\fA R) = \beta_i(\fA)$.
\end{example}

\begin{example}[The boundary condition cannot simply be omitted for complete domains]\label{domain}
Let $B := K[[x,y,z]]$, and embed $B$ in
$D := K[[x,y,u]]$ by sending $z$ to $u^2$. Put $\fb=(x,y)B$ and
$R := B\oplus \fb u\subseteq D = B\oplus Bu$.
This is a subring, since $(\fb u)^2=z\fb^2\subseteq B$. It is a module-finite
$B$-algebra and a local domain with maximal ideal 
$\m_R = (x,y,z,xu,yu)R$ and  residue field $K$. Indeed, 
an element $a+bu$, with $a\in B$
and $b\in \fb $, is a unit whenever $a$ is a unit, since
$\displaystyle (a+bu)^{-1}=\frac{a-bu}{a^2-b^2z}\in R$.
The ring $R$ is complete, because it is module-finite over the complete local
ring $B$ and $I:=(x,y,z)R$ is $\m_R$-primary. Thus $R$ is a
complete local domain of dimension $3$, and $x,y,z$ is a system of
parameters. Its parameter multiplicity is
$\e_I(R)=\rank_B(R)=2$. Here the first equality is the usual rank formula for multiplicity over the regular local ring $B$.  Note again that the residue fields agree, so $\e_I(R)$ agrees with the multiplicity of $R$ as a $B$-module with respect to the maximal ideal of $B$.

For $t \geq 2$, consider the formal monomial ideal
\[
 \fA_t=(X^t,Y^t,XZ^t,YZ^t,Z^{t^5})
 \inc\Z[X,Y,Z].
\]
$\fB(\fA_t)$ consists of $Z^t,Z^{t+1},\ldots,Z^{t^5-1}$ along
with the monomials $X^iY^jZ^k$ that correspond to the cube
$\{(i,j,k) \in \N^3: 0\leq i,\, j,\, k < t\}.$ 
 Consequently, $\cN(\fA_t)=\la_B(B/\fA_t B) = t^3+t^5-t$.  
 Note that $\cN_0(\fA_t) = t^5+3t^2-4t+1$ and thus 
 $\gamma(\fA_t)\to 1$ as $t \to \infty$, so that this family does \textit{not} 
have asymptotically sparse boundary.  

Using the direct sum decomposition of $R$, we obtain that 
$\la(R/\fA_tR)
   =\la_B(B/\fA_t B)+\la_B(\fb /\fA_t \fb )$.
To compute the second term, let $\theta: B \to B^2$  be the
$B$-linear map such that $1 \mapsto (-y,\, x)$ and apply $\blank \otimes_B B/\fA_t B$ 
to the presentation $\displaystyle B\arwf{\theta}B^2 \to \fb \to 0$.
The resulting first map has kernel
$0:_{B/\fA_tB}(x,y)$, so that 
\[
 \la_B(\fb /\fA_t \fb )
   =\cN(\fA_t)+\la_B(0:_{B/\fA_tB}(x,y)).
\]
A monomial basis for this annihilator is
\[
 \{z^c:t \leq c < t^5\}
 \ \cup\
 \{x^{t-1}y^{t-1}z^c : 0\leq c < t\}.
\]
It has $(t^5 - t) + t = t^5$ elements. Therefore the exact length formula is
\[
 \la(R/\fA_tR)  = \la_B(B/\fA_tB) + \la_B(\fb/\fA_t\fb)= 2\cN(\fA_t) + t^5,
\]
and hence 
$\displaystyle \lim_{t\to\infty} \frac{\la(R/\fA_tR)}{\cN(\fA_t)}  =3
 \neq 2 = \e_I(R)$. Observe as well that $\fA_{t+1}\inc\fA_t$ \linebreak \vspace*{3pt}
and $\fA_t\subseteq(X,Y,Z)^t$.  Consequently, the ideals $\fA_t R$ form a 
descending sequence cofinal with the powers of  $\m$.  
Each is $\m$-primary, and $\fA_t R \inc I^t \inc \m^t.$
\end{example}

\begin{remark} In Example~\ref{domain} we have used $t^5$ for definiteness, but one
obtains a similar example using instead $t^k$ for every integer $k \geq 3$.   
If one makes this replacement for an integer $k\geq4$,
$\displaystyle
\frac{\la_R(R/\fA_tR)}{\cN(\fA_t)}
=2+\frac{1}{1+t^{3-k}-t^{1-k}} \to 3$, while $\gamma(\fA_t) \to 1.$
If one replaces $t^5$ by $t^3$, the length ratio instead tends to $5/2$, while
$\gamma(\fA_t) \to 1/2$.
Hence, for all $k \geq 3$, the length limits differ from $\e_I(R)=2$.
For each fixed $k\geq 3$, one still has $\fA_t \inc (X,Y,Z)^t$, and their extensions
to $R$ are $\m$-primary, so that their extensions remain cofinal
with the powers of $\m$.
\end{remark}

\quad\bigskip

\quad\bigskip

$\begin{array}{ll}
\hspace*{-18pt} \textrm{Department of Mathematics}&\qquad\quad \textrm{Department of Mathematics and Statistics}\\
\hspace*{-18pt}\textrm{University of Michigan}   &\qquad\quad\textrm{University of New Mexico}\\
\hspace*{-18pt}\textrm{Ann Arbor, MI 48109--1043}      &\quad\qquad \textrm{Albuquerque, New Mexico 87131}\\
\hspace*{-18pt}\textrm{USA}                            &\quad\qquad \textrm{USA}\\
\hspace*{-18pt}\quad & \quad\\
\hspace*{-18pt}\textrm{E-mail: hochster@umich.edu}     &\qquad\quad \textrm{E-mail: ostrahan@unm.edu}\\ 
 
\end{array}$

\end{document}